\documentclass{article}

\input{diagrams}
\usepackage{epsfig}

\newcommand{\mcm}[3]{\newcommand{#1}[#2]{{\ensuremath{#3}}}}

\usepackage{latexsym}
\usepackage{amssymb}

\mcm{\diagspace}{0}{\mbox{\hspace{2em}}}

\mcm{\cat}{1}{\mc{#1}}
\mcm{\fcat}{1}{\mb{#1}}
\mcm{\mc}{1}{\mathcal{#1}}
\mcm{\mr}{1}{\mathrm{#1}}
\mcm{\mb}{1}{\mathbf{#1}}
\newcommand{\url}[1]{\mbox{\tt #1 }}

\mcm{\of}{0}{\raisebox{0.2mm}{\ensuremath{\scriptstyle\circ}}}
\mcm{\Hom}{0}{\mr{Hom}}

\mcm{\go}{0}{\rTo}
\mcm{\goby}{1}{\rTo^{#1}}

\newcommand{\piccy}[1]{\epsfig{file=#1}}

\mcm{\Top}{0}{\fcat{Top}}

\newcommand{\demph}[1]{\emph{#1}}

\newcommand{\ditto}{,,}
\mcm{\blob}{0}{\raisebox{.3ex}{\ensuremath{\scriptscriptstyle{\bullet}}}}
\newcommand{\place}[3]{\put(#1,#2){\makebox(0,0)[c]{#3}}}

\newlength{\gwidth}	
\newlength{\gvert}	
\newlength{\gdrop}	
\newlength{\gbaredrop}	
\newlength{\goffset}	
\newlength{\gtemp}	

\newcommand{\present}[1]{%
\makebox[1\gwidth]{%
\rule[-1\gdrop]{0ex}{1\gvert}%
\raisebox{-1\gbaredrop}{#1}}}

\newcommand{\presentl}[1]{%
\makebox[1\gwidth][l]{%
\rule[-1\gdrop]{0ex}{1\gvert}%
\raisebox{-1\gbaredrop}{#1}}}

\newcommand{\presentr}[1]{%
\makebox[1\gwidth][r]{%
\rule[-1\gdrop]{0ex}{1\gvert}%
\raisebox{-1\gbaredrop}{#1}}}

\newcommand{\ginitdims}[2]{
\setlength{\unitlength}{1em}
\setlength{\goffset}{.25\unitlength}
\setlength{\gwidth}{#1\unitlength}
\setlength{\gvert}{#2\unitlength}
\setlength{\gdrop}{.5\gvert}
\addtolength{\gdrop}{-1\goffset}
\setlength{\gbaredrop}{1\gdrop}
\addtolength{\gvert}{.6\unitlength}
\addtolength{\gdrop}{.3\unitlength}}	

\newcommand{\cinitdims}[2]{
\setlength{\unitlength}{1em}
\setlength{\goffset}{.35\unitlength}
\setlength{\gwidth}{#1\unitlength}
\setlength{\gvert}{#2\unitlength}
\setlength{\gdrop}{.5\gvert}
\addtolength{\gdrop}{-1\goffset}
\setlength{\gbaredrop}{1\gdrop}
\addtolength{\gvert}{.6\unitlength}
\addtolength{\gdrop}{.3\unitlength}}	

\newcommand{\gsinitdims}[2]{
\setlength{\unitlength}{0.5em}
\setlength{\goffset}{.25\unitlength}
\setlength{\gwidth}{#1\unitlength}
\setlength{\gvert}{#2\unitlength}
\setlength{\gdrop}{.5\gvert}
\addtolength{\gdrop}{-1\goffset}
\setlength{\gbaredrop}{1\gdrop}
\addtolength{\gvert}{.6\unitlength}
\addtolength{\gdrop}{.3\unitlength}}	

\newcommand{\sidespic}[1]{%
\settowidth{\gtemp}{\ensuremath{#1}}%
\addtolength{\gwidth}{1\gtemp}}

\newcommand{\abovepic}[1]{%
\settoheight{\gtemp}{\ensuremath{#1}}%
\addtolength{\gvert}{1\gtemp}%
\settodepth{\gtemp}{\ensuremath{#1}}%
\addtolength{\gvert}{1\gtemp}}

\newcommand{\belowpic}[1]{%
\settoheight{\gtemp}{\ensuremath{#1}}%
\addtolength{\gvert}{1\gtemp}%
\addtolength{\gdrop}{1\gtemp}%
\settodepth{\gtemp}{\ensuremath{#1}}%
\addtolength{\gvert}{1\gtemp}%
\addtolength{\gdrop}{1\gtemp}}

\newcommand{\cell}[4]{\put(#1,#2){\makebox(0,0)[#3]{\ensuremath{#4}}}}
\mcm{\zmark}{0}{\scriptstyle{\bullet}}

\newcommand{\pregfst}[1]{%
\begin{picture}(0.5,0.2)(-0.5,-0.2)%
\cell{-0.1}{-0.2}{tr}{#1}%
\cell{0}{0}{c}{\zmark}%
\end{picture}}

\mcm{\gfst}{1}{%
\ginitdims{0.5}{0.4}%
\sidespic{#1}%
\belowpic{#1}%
\presentr{\pregfst{#1}}}

\newcommand{\preglst}[1]{%
\begin{picture}(0.5,0.2)(0,-0.2)%
\cell{0.1}{-0.2}{tl}{#1}%
\cell{0.05}{0}{c}{\zmark}%
\end{picture}}

\mcm{\glst}{1}{%
\ginitdims{.5}{.4}%
\sidespic{#1}%
\belowpic{#1}%
\presentl{\preglst{#1}}}

\newcommand{\preglft}[1]{%
\begin{picture}(0,0.2)(0,-0.2)%
\cell{-0.1}{-0.2}{tr}{#1}%
\cell{0.05}{0}{c}{\zmark}%
\end{picture}}

\mcm{\glft}{1}{%
\ginitdims{0}{.4}%
\belowpic{#1}%
\present{\preglft{#1}}}

\newcommand{\pregrgt}[1]{%
\begin{picture}(0,0.2)(0,-0.2)%
\cell{0.1}{-0.2}{tl}{#1}%
\cell{0.05}{0}{c}{\zmark}%
\end{picture}}

\mcm{\grgt}{1}{%
\ginitdims{0}{.4}%
\belowpic{#1}%
\present{\pregrgt{#1}}}

\newcommand{\pregblw}[1]{%
\begin{picture}(0,0.3)(0,-0.3)
\cell{0}{-0.3}{t}{#1}%
\cell{0.05}{0}{c}{\zmark}%
\end{picture}}

\mcm{\gblw}{1}{%
\ginitdims{0}{.6}%
\belowpic{#1}%
\present{\pregblw{#1}}}

\newcommand{\pregfbw}[1]{%
\begin{picture}(0,0.65)(0,-0.65)
\cell{0}{-0.65}{t}{#1}%
\cell{0.05}{0}{c}{\zmark}%
\end{picture}}

\mcm{\gfbw}{1}{%
\ginitdims{0}{1.3}%
\belowpic{#1}%
\present{\pregfbw{#1}}}

\newcommand{\pregzero}[1]{%
\begin{picture}(0.8,0.4)(-0.4,-0.4)
\cell{0}{-0.4}{t}{#1}%
\cell{0}{0}{c}{\zmark}%
\end{picture}}

\mcm{\gzero}{1}{%
\ginitdims{0.8}{.6}%
\belowpic{#1}%
\sidespic{#1}%
\present{\pregzero{#1}}}

\newcommand{\pregone}[1]{%
\begin{picture}(5,0.4)(0,-0.2)%
\cell{2.5}{0.2}{b}{#1}%
\put(0,0){\vector(1,0){5}}%
\end{picture}}

\mcm{\gone}{1}{%
\ginitdims{5}{0.4}%
\abovepic{#1}%
\present{\pregone{#1}}}

\newcommand{\pregtwo}[3]{%
\begin{picture}(5,3.4)(0,-0.2)%
\cell{2.5}{3.2}{b}{#1}%
\cell{2.5}{-.2}{t}{#2}%
\cell{2.7}{1.5}{l}{#3}%
\qbezier(0,1.5)(2.5,4.5)(5,1.5)%
\qbezier(0,1.5)(2.5,-1.5)(5,1.5)%
\put(5,1.5){\vector(1,-1){0}}%
\put(5,1.5){\vector(1,1){0}}%
\put(2.5,2.5){\vector(0,-1){2}}%
\end{picture}}

\mcm{\gtwo}{3}{%
\ginitdims{5}{3.4}%
\abovepic{#1}%
\belowpic{#2}%
\present{\pregtwo{#1}{#2}{#3}}}

\newcommand{\pregthree}[5]{%
\begin{picture}(5,5.4)(0,-1.2)%
\cell{2.5}{4.2}{b}{#1}%
\cell{1.5}{1.7}{b}{#2}%
\cell{2.5}{-1.2}{t}{#3}%
\cell{2.7}{2.75}{l}{#4}%
\cell{2.7}{0.25}{l}{#5}%
\qbezier(0,1.5)(2.5,6.5)(5,1.5)%
\qbezier(0,1.5)(2.5,-3.5)(5,1.5)%
\put(0,1.5){\vector(1,0){5}}%
\put(2.5,3.5){\vector(0,-1){1.5}}%
\put(2.5,1){\vector(0,-1){1.5}}%
\put(5,1.5){\vector(1,-3){0}}%
\put(5,1.5){\vector(1,3){0}}%
\end{picture}}

\mcm{\gthree}{5}{%
\ginitdims{5}{5.4}%
\abovepic{#1}%
\belowpic{#3}%
\present{\pregthree{#1}{#2}{#3}{#4}{#5}}}

\mcm{\gzersu}{0}{%
\gsinitdims{0}{.6}%
\present{\pregblw{}}}

\mcm{\gonesu}{0}{%
\gsinitdims{5}{0.4}%
\present{\pregone{}}}

\mcm{\gtwosu}{0}{%
\gsinitdims{5}{3.4}%
\present{\pregtwo{}{}{}}}

\mcm{\gthreesu}{0}{%
\gsinitdims{5}{5.4}%
\present{\pregthree{}{}{}{}{}}}

\newcommand{\precone}[1]{%
\begin{picture}(4.2,0.4)(-0.3,-0.2)%
\cell{1.8}{0.2}{b}{#1}%
\put(0,0){\vector(1,0){3.6}}%
\end{picture}}

\mcm{\cone}{1}{%
\cinitdims{4.2}{0.4}%
\abovepic{#1}%
\present{\precone{#1}}}

\mcm{\gfstsu}{0}{%
\gsinitdims{0.5}{0.4}%
\presentr{\pregfst{}}}

\mcm{\glstsu}{0}{%
\gsinitdims{0.5}{0.4}%
\presentl{\preglst{}}}

\newcommand{\prectwodbl}[3]%
{\begin{picture}(4.2,3.4)(-0.1,-0.2)%
\cell{2}{3.2}{b}{#1}%
\cell{2}{-0.2}{t}{#2}%
\cell{2.3}{1.5}{l}{#3}%
\qbezier(0,2)(2,4)(4,2)%
\qbezier(0,1)(2,-1)(4,1)%
\put(4,2){\vector(1,-1){0}}%
\put(4,1){\vector(1,1){0}}%
\put(1.9,2.5){\line(0,-1){1.8}}%
\put(2.1,2.5){\line(0,-1){1.8}}%
\cell{2.01}{0.4}{b}{\vee}%
\end{picture}}

\mcm{\ctwodbl}{3}{%
\cinitdims{4.2}{3.4}%
\abovepic{#1}%
\belowpic{#2}%
\present{\prectwodbl{#1}{#2}{#3}}}

\newcommand{\precthreedbl}[5]{%
\begin{picture}(4.2,5.4)(-0.1,-0.2)%
\cell{2}{5.2}{b}{#1}%
\cell{1}{2.7}{b}{#2}%
\cell{2}{-.2}{t}{#3}%
\cell{2.3}{3.75}{l}{#4}%
\cell{2.3}{1.25}{l}{#5}%
\qbezier(0,3)(2,7)(4,3)%
\qbezier(0,2)(2,-2)(4,2)%
\put(0,2.5){\vector(1,0){4}}%
\put(1.9,4.5){\line(0,-1){1.3}}%
\put(2.1,4.5){\line(0,-1){1.3}}%
\cell{2.01}{2.9}{b}{\vee}%
\put(1.9,2){\line(0,-1){1.3}}%
\put(2.1,2){\line(0,-1){1.3}}%
\cell{2.01}{0.4}{b}{\vee}%
\put(4,3){\vector(1,-3){0}}%
\put(4,2){\vector(1,3){0}}%
\end{picture}}

\mcm{\cthreedbl}{5}{%
\cinitdims{4.2}{5.4}%
\abovepic{#1}%
\belowpic{#3}%
\present{\precthreedbl{#1}{#2}{#3}{#4}{#5}}}

\newcommand{\precthreecelltrp}[5]{%
\begin{picture}(8.2,5)(-4.1,-2.5)%
\cell{0}{2.5}{b}{#1}%
\cell{0}{-2.5}{t}{#2}%
\cell{-1.8}{0}{r}{#3}%
\cell{1.8}{0}{l}{#4}%
\cell{0}{0.3}{b}{#5}%
\qbezier(-4,0.5)(0,4)(4,0.5)%
\qbezier(-4,-0.5)(0,-4)(4,-0.5)%
\qbezier(-0.6,2)(-2.6,0)(-0.6,-2)%
\qbezier(-0.4,2)(-2.4,0)(-0.5,-1.9)%
\cell{-0.6}{-2}{b}{\lrcorner}%
\qbezier(0.4,2)(2.4,0)(0.5,-1.9)%
\qbezier(0.6,2)(2.6,0)(0.6,-2)%
\cell{0.65}{-2}{b}{\llcorner}%
\put(-1,0.15){\line(1,0){1.7}}%
\put(-1,0){\line(1,0){2}}%
\put(-1,-0.15){\line(1,0){1.7}}%
\cell{1.15}{0}{r}{>}%
\put(4,0.5){\vector(1,-1){0}}%
\put(4,-0.5){\vector(1,1){0}}%
\end{picture}}

\mcm{\cthreecelltrp}{5}{%
\cinitdims{8.2}{5}%
\abovepic{#1}%
\belowpic{#2}%
\present{\precthreecelltrp{#1}{#2}{#3}{#4}{#5}}}

\newcommand{\pregtwodblsu}{%
\begin{picture}(5,3.4)(0,-0.2)%
\qbezier(0,1.5)(2.5,4.5)(5,1.5)%
\qbezier(0,1.5)(2.5,-1.5)(5,1.5)%
\put(5,1.5){\vector(1,-1){0}}%
\put(5,1.5){\vector(1,1){0}}%
\cell{2.5}{1.5}{c}{\Downarrow}%
\end{picture}}

\mcm{\gtwodblsu}{0}{%
\gsinitdims{5}{3.4}%
\present{\pregtwodblsu}}

\newcommand{\pregthreedblsu}{%
\begin{picture}(5,5.4)(0,-1.2)%
\qbezier(0,1.5)(2.5,6.5)(5,1.5)%
\qbezier(0,1.5)(2.5,-3.5)(5,1.5)%
\put(0,1.5){\vector(1,0){5}}%
\cell{2.5}{2.75}{c}{\Downarrow}%
\cell{2.5}{0.25}{c}{\Downarrow}%
\put(5,1.5){\vector(1,-3){0}}%
\put(5,1.5){\vector(1,3){0}}%
\end{picture}}

\mcm{\gthreedblsu}{0}{%
\gsinitdims{5}{5.4}%
\present{\pregthreedblsu}}

\begin{document}

\title{Topology and Higher-Dimensional Category Theory: the Rough Idea}
\author{Tom Leinster}
\date{}

\maketitle

\vspace{-5mm}

\begin{center}	\bfseries
Abstract
\end{center}
\begin{quotation}
Higher-dimensional category theory is the study of $n$-categories, operads,
braided monoidal categories, and other such exotic structures.  Although it
can be treated purely as an algebraic subject, it is inherently topological
in nature: the higher-dimensional diagrams one draws to represent these
structures can be taken quite literally as pieces of topology.  Examples of
this are the braids in a braided monoidal category, and the pentagon which
appears in the definitions of both monoidal category and $A_\infty$-space.

I will try to give a Friday-afternoonish description of some of the dreams
people have for higher-dimensional category theory and its interactions with
topology.  Grothendieck, for instance, suggested that tame topology should be
the study of $n$-groupoids; others have hoped that an $n$-category of
cobordisms between cobordisms between \ldots\ will provide a clean setting for
TQFT; and there is convincing evidence that the whole world of $n$-categories
is a mirror of the world of homotopy groups of spheres.
\end{quotation}

\vspace{3mm}

\emph{%
These are notes from talks given in London and Sussex in summer 2001.  I
thank Dicky Thomas and Roger Fenn for their invitations and the audiences
for their comments, many of which are incorporated here.}

\vspace{5mm}

What I want to give you in the next hour is an informal description of what
higher-dimensional category theory is and might be, and how it is relevant to
topology.  There will be no real theorems, proofs or definitions.  But to
whet your appetite, here's a question which we'll reach an answer to by the
end:

\paragraph*{Question} What is the close connection between the following two
facts? 
\begin{description}
\item[A] No-one ever got into trouble for leaving out the brackets in a
tensor product of several objects (abelian groups, etc.).  For instance, it's
safe to write $A\otimes B\otimes C$ instead of $(A\otimes B)\otimes C$ or
$A\otimes (B\otimes C)$.
\item[B] There exist non-trivial knots (that is, knots which cannot be
undone) in $\Bbb{R}^3$.
\end{description}

\section{The Very Rough Idea}

In ordinary category theory we have diagrams of objects and arrows such as
\[
\gfstsu\ \gonesu\ \gzersu\ \gonesu\ \gzersu\ \gonesu\ \glstsu.
\]
We can also consider more complex category-like structures, in which there
are diagrams such as
\[
\piccy{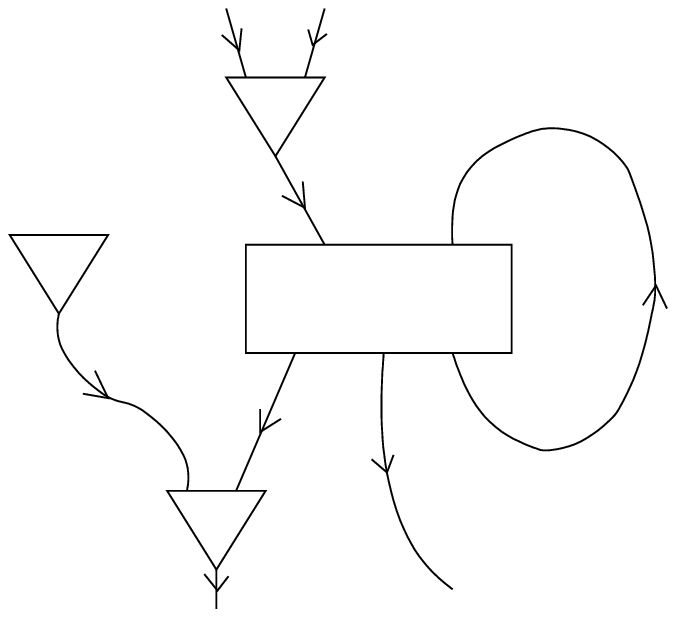}.
\]
This looks like an electronic circuit diagram or a flow chart; the unifying
idea is that of `information flow'. It can be redrawn as
\[
\piccy{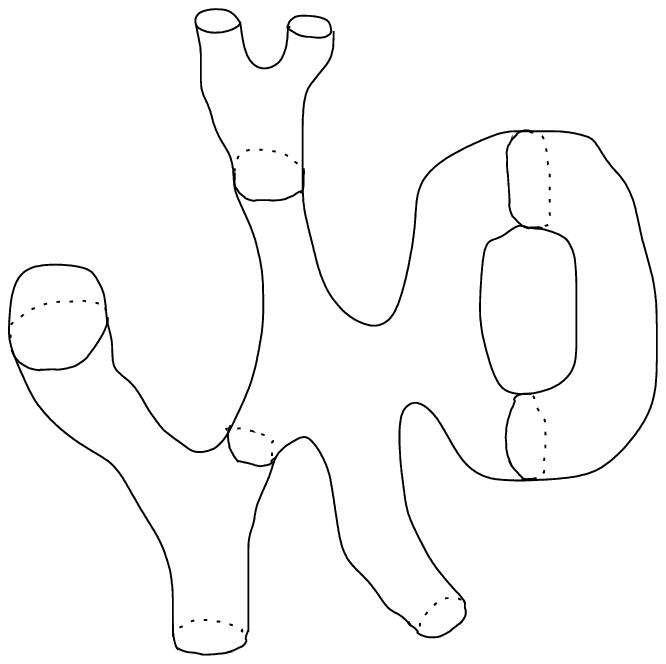},
\]
which looks like a surface or a diagram from topological quantum field
theory.

You can also use diagrams like this to express algebraic laws, e.g.\
commutativity: 
\[
\piccy{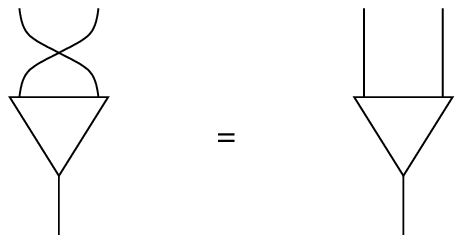}.
\]
This is perhaps more clear when labels are added:
\[
\piccy{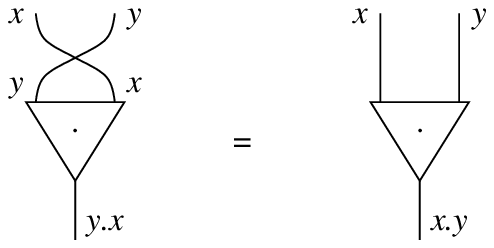}.
\]
The fact that two-dimensional TQFTs are essentially Frobenius algebras
is an example of an explicit link between the spatial and algebraic
aspects of diagrams like these.

Moreover, if you allow crossings, as in the diagrams for commutativity or as
in
\[
\piccy{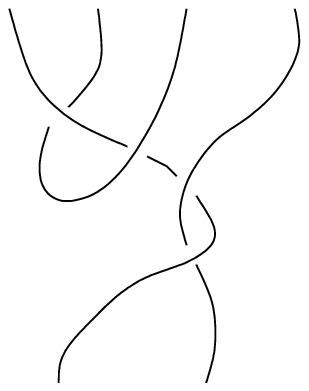},
\]
then you start getting pictures that look like knots; and there are indeed
well-established relations between knot theory and higher categorical
structures.

So the idea is:
\begin{trivlist} \item
\fbox{\parbox{0.98\textwidth}{\centering%
in ordinary category theory we have $1$-dimensional arrows $\go$;\\
in higher-dimensional category theory we have higher-dimensional arrows.}}
\end{trivlist}
The natural geometry of these higher-dimensional arrows is what makes
higher-dimensional category theory an inherently topological subject.

\paragraph*{Example} For a concrete example of what I mean by a
`higher-dimensional arrow', consider operads.  An operad consists of
operations shaped like
\[
\piccy{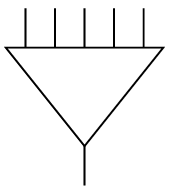},
\]
with many inputs (five, here) and one output.  Any tree of operations, such
as
\[
\piccy{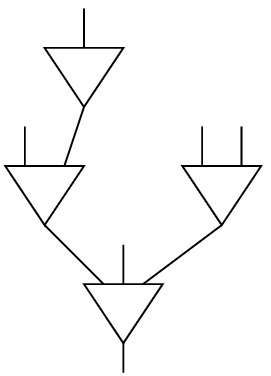},
\]
can be composed to give a single operation.  (Strictly speaking, these
pictures are appropriate for operads without a symmetric group action: `non-$\Sigma$ operads', or `planar operads' as they are sometimes called.)

\paragraph*{} There's a whole zoo of structures coming under the heading of
higher-dimensional category theory, including operads, generalized operads
(of which the variety familiar to topologists is a basic special case),
multicategories, various flavours of monoidal categories, and $n$-categories.
Today I'll concentrate on $n$-categories.

Terminology: a `higher-dimensional category' or $n$-category is not a
special kind of category, but a generalization of the notion of category;
compare the usage of `quantum group'.  A $1$-category is the same thing
as an ordinary category.

\section{$n$-Categories}

Here's a very informal

\paragraph*{`Definition'} Let $n\geq 0$.  An
\demph{$n$-category} consists of 
\begin{itemize}
\item \demph{0-cells} or \demph{objects}, $A, B, \ldots$
\item \demph{1-cells} or \demph{morphisms}, drawn as
$A \cone{f} B$ 
\item \demph{2-cells} $A \ctwodbl{f}{g}{\alpha} B$ (`morphisms between
morphisms') 
\item \demph{3-cells}
$A \cthreecelltrp{f}{g}{\alpha}{\beta}{\Gamma} B$ 
(where the arrow
labelled $\Gamma$ is meant to be going in a direction perpendicular to the
plane of the paper)
\item \ldots
\item all the way up to \demph{$n$-cells}
\item various kinds of \demph{composition}, e.g.
\[
\begin{array}{ccc}
A \cone{f} B \cone{g} C					&\textrm{gives}
&A \cone{g\of f} C			\\
&&\textrm{(as we're all familiar with)}	\\
A \cthreedbl{f}{g}{h}{\alpha}{\beta} B			&\textrm{gives}
&A \ctwodbl{f}{h}{\beta\of\alpha} B	\\
A \ctwodbl{f}{g}{\alpha} A' \ctwodbl{f'}{g'}{\alpha'}
A'' 							&\textrm{gives}
&A \ctwodbl{f'\of f}{g'\of g}{\!\!\!\!\!\!\!\!\alpha' *\alpha} A''\\
&&\textrm{(the $*$ notation is traditional}\\
&&\textrm{but not particularly well-chosen)}
\end{array}
\]
and so on in higher dimensions (which I won't attempt to draw); and similarly
\demph{identities}.  
\end{itemize}
These compositions are required to `all fit together nicely'---a phrase
hiding many subtleties. 

\demph{$\infty$-categories} (also known as \demph{$\omega$-categories}) are
defined similarly, by going on up the dimensions forever instead of stopping
at $n$.

\paragraph*{} You could also consider some kind of higher categorical
structure which involved cubical (or other) shapes, e.g. a 2-cell might look
like
\begin{diagram}[size=2em,abut]
\bullet	&\rTo		&\bullet\\
\dTo	&\Downarrow	&\dTo	\\
\bullet	&\rTo		&\bullet.	
\end{diagram}
This is one of the `zoo of structures' mentioned earlier, but is something
other than an $n$-category.

\paragraph*{Critical Example} This is the excuse for putting the word
`Grothendieck' into the abstract.  Any topological space $X$ gives rise to an
$\infty$-category $\Pi_\infty(X)$ (its \demph{fundamental
$\infty$-groupoid}), in which
\begin{itemize}
\item 0-cells are points of $X$, drawn as $\ \blob$
\item 1-cells are paths in $X$ (parametrized, i.e. maps $[0,1] \go X$), drawn
as $\gfstsu\gonesu\glstsu$ ---though whether that's meant to be a picture in
the space $X$ or the $\infty$-category $\Pi_\infty(X)$ is deliberately
ambiguous; we're trying to blur the distinction between geometry and algebra
\item 2-cells are homotopies of paths (relative to endpoints), drawn
$\gfstsu\gtwodblsu\glstsu$ 
\item 3-cells are homotopies of homotopies of paths (i.e.\ suitable maps
$[0,1]^3 \go X$) 
\item \ldots
\item composition is by pasting paths and homotopies.
\end{itemize}
So $\Pi_\infty(X)$ should contain all the information you want about $X$ if
your context is `tame topology'; e.g.\ you should be able to compute from it
the homotopy, homology and cohomology of $X$.  (The word `groupoid' means
that all cells of dimension $>0$ are invertible.)

You can also truncate after $n$ steps in order to obtain $\Pi_n(X)$, the
\demph{fundamental $n$-groupoid} of $X$: e.g.\ $\Pi_1(X)$ is the familiar
fundamental groupoid.

\paragraph{Alert} As you may have noticed, composition in $\Pi_\infty(X)$
isn't genuinely associative; nor is it unital, and nor are the cells
genuinely invertible (only up to homotopy).  We're therefore interested in
\demph{weak} $n$-categories, where the `fitting together nicely' only happens
up to some kind of equivalence, rather than \demph{strict} $n$-categories,
where associativity etc.\ hold in the strict sense. 

To define strict $n$-categories precisely turns out to be easy.  To define
weak $n$-categories, we face the same kind of challenge as algebraic
topologists did in the 60s, when they were trying to state the exact sense in
which a loop space is a topological group.  It clearly isn't a group in the
literal sense, as composition of paths isn't associative; but it is
associative up to homotopy, and if you pick specific homotopies to do this
job then these homotopies obey laws of their own---or at least, obey them up
to homotopy; and so on.  At least two precise formulations of `group up to
(higher) homotopy' became popular: Stasheff's $A_\infty$-spaces and Segal's
special $\Delta$-spaces.  (More exactly, these are notions of monoid or
semigroup up to homotopy; the inverses are dealt with separately.)  

The situation for weak $n$-categories is similar but more extreme: there are
something like a dozen proposed definitions that I know of, and no-one has
much idea of how they relate to one another.  (Some of the reasons for this
chaos are good: there are real conceptual difficulties in saying what it
means for two definitions of weak $n$-category to be equivalent.)  Happily,
we can ignore all this today and work informally.  This means that nothing I
say from now on is true with any degree of certainty or accuracy.

At this point you might be thinking: can't we do away with this difficult
theory of weak $n$-categories and just stick to the strict ones?  The answer
is: if you're interested in topology, no.  The difference between the weak
and strict theories is genuine and nontrivial: for while it is true that
every weak $2$-category is equivalent to some strict one, and so it is also
true that homotopy $2$-types can be modelled by strict $2$-groupoids, neither
of these things is true in dimensions $\geq 3$.  For instance, there exist
spaces $X$ such that the weak $3$-category $\Pi_3(X)$ is not equivalent to
any strict $3$-category.

From now on, `$n$-category' will mean `weak $n$-category'.  The strict ones
hardly arise in nature.

\paragraph*{Some More Examples} of $\infty$-categories:
\begin{description}
\item[\Top] This is very similar to the $\Pi_\infty$ example above.  \Top\
has:
\begin{itemize}
\item 0-cells: topological spaces
\item 1-cells: continuous maps
\item 2-cells $X \ctwodbl{f}{g}{} Y$: homotopies between $f$ and $g$
\item 3-cells: homotopies between homotopies (i.e.\ suitable maps $[0,1]^2
\times X \go Y$) 
\item \ldots
\item composition as expected.
\end{itemize}

\item[\fcat{ChCx}] This $\infty$-category has:
\begin{itemize}
\item 0-cells: chain complexes (of abelian groups, say)
\item 1-cells: chain maps
\item 2-cells: chain homotopies
\item 3-cells $A\cthreecelltrp{f}{g}{\alpha}{\beta}{\Gamma}B$: homotopies
between homotopies, i.e.\ maps $\Gamma: A \go B$ of degree $2$ such that
$d\Gamma - \Gamma d = \beta - \alpha$
\item \ldots
\item composition: more or less as expected, but some choices are involved.
For instance, if you try to write down the composite of two chain homotopies
$\gfstsu\gtwodblsu\gzersu\gtwodblsu\glstsu$ then you'll find that there are two
equally reasonable ways of doing it: one `left-handed', one `right-handed'.
This is something like choosing the parametrization when deciding how to
compose two loops in a space (usual choice: do everything at double speed).
Somehow the fact that there's no canonical choice means that weakness of the
resulting $\infty$-category is inevitable. 
\end{itemize}
In a reasonable world there ought to be some kind of map $\fcat{Chains}:
\fcat{Top} \go \fcat{ChCx}$.

\item[Bord] This is an $\infty$-category of (co?)bordisms.  
\begin{itemize}
\item 0-cells: 0-manifolds, where `manifold' means `compact, smooth, oriented
manifold'.  A typical 0-cell is $\ \blob\ \ \blob\ \ \blob\ \ \blob\ $.
\item 1-cells: 1-manifolds with corners, i.e. cobordisms between 0-manifolds,
such as
\[
\piccy{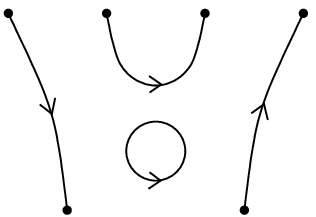}
\]
(this being a 1-cell from the 4-point 0-manifold to the 2-point 0-manifold).
In Atiyah-Segal-style TQFT, we'd stop here and take \emph{isomorphism
classes} of the 1-cells just described, to make a category.  We avoid this
(unnatural?) quotienting out and carry on up the dimensions.
\item 2-cells: 2-manifolds with corners, such as
\[
\begin{array}[c]{c}
\piccy{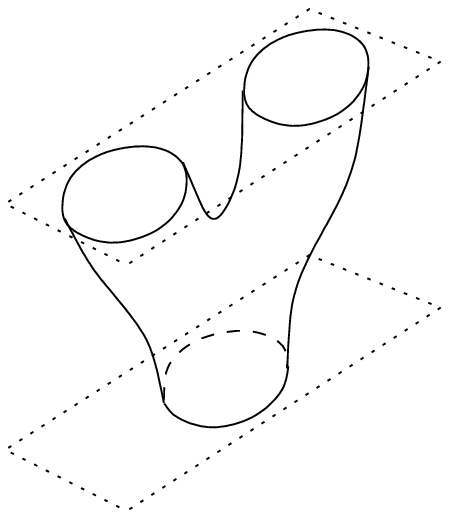}	
\end{array}
\ \ \textrm{or}\ \ 	
\begin{array}[c]{c}
\piccy{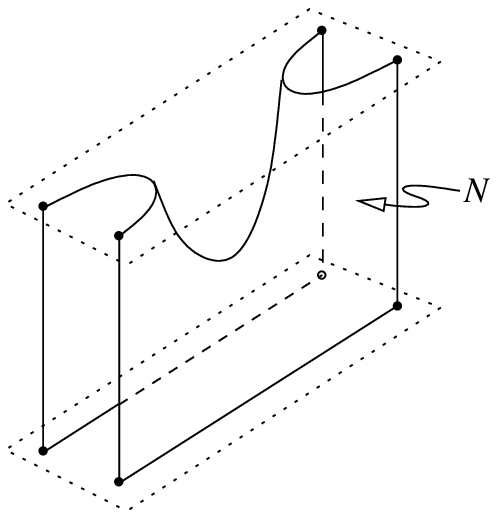}	
\end{array}
\]
Here $N$ is a 2-cell $L \ctwodbl{M}{M'}{N} L'$, where 
\[
L=L'=\ \blob\ \ \blob\ , \diagspace
M=\begin{array}[c]{c}
\piccy{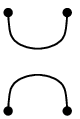}
\end{array},
\diagspace
M'=\begin{array}[c]{c}
\piccy{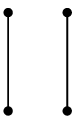}
\end{array}.  
\]
I've left the orientations off the pictures, but $N$ is meant to be oriented
so as to agree with the orientations of $M$ and $M'$.  In his~\cite{Khov}
paper, Khovanov discusses TQFTs with corners in the language of 2-categories;
he'd stop here and take isomorphism classes of the 2-cells just described, to
make a 2-category.  Again, we do not quotient out but keep going up the
dimensions.
\item 3-cells, 4-cells, \ldots\ are defined similarly
\item composition is gluing of manifolds.
\end{itemize}
Some authors discuss `extended TQFTs' via the notion of $n$-vector space.
A $0$-vector space is a complex number; a $1$-vector space is an ordinary
vector space; and $n$-vector spaces for higher $n$ are something more
sophisticated.  I won't attempt to say anything about this; see Further
Reading below for pointers.  
\end{description}

\section{Degenerate $n$-Categories}

So far we've seen that topological structures provide various good examples
of $n$-categories, and that alone might be enough to convince you that
$n$-categories are interesting from a topological point of view.  But the
relationship between topology and higher-dimensional category theory is
actually much more intimate than that.  To see this, we'll consider
$n$-categories which are degenerate in various ways.  This doesn't sound very
promising, and it'll seem at first as if it's a purely formal exercise, but
in a little while the intrinsic topology should begin to shine through.

\paragraph*{Some Degeneracies}

\begin{itemize}

\item A category \cat{C} with only one object is the same thing as a monoid
($=$ semigroup with unit) $M$.
For if the single object of \cat{C} is called $*$, say, then \cat{C} just
consists of the set $\Hom(*,*)$ together with a binary operation of
\marginpar{\hspace{1em}\piccy{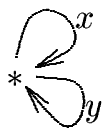}}
composition and a unit element $1$, obeying the usual axioms.  So we have:
\begin{eqnarray*}
\textrm{morphism in }\cat{C}	&=	&\textrm{element of }M		\\
\of \textrm{ in } \cat{C}	&=	&\cdot \textrm{ in }M.
\end{eqnarray*}

\pagebreak

\item A 2-category \cat{C} with only one $0$-cell is the same thing as a
monoidal category \cat{M}.  (Private thought: if \cat{C} has only one
$0$-cell then there are only interesting things happening in the top two
dimensions, so it's going to be \emph{some} kind of one-dimensional
structure.)  This works as follows:
\begin{eqnarray*}
1\textrm{-cell in }\cat{C}	&=	&\textrm{object of }\cat{M}	\\
2\textrm{-cell in }\cat{C}	&=	&\textrm{morphism of }\cat{M}	\\
\textrm{composition } \ \gzersu\gonesu\gzersu\gonesu\gzersu\  \textrm{ in }
\cat{C} 
&=	
&\otimes \textrm{ of objects in }\cat{M}				\\
\textrm{composition } \ \gzersu\gthreedblsu\gzersu\  \textrm{ in } \cat{C}
&=	
&\of \textrm{ of morphisms in }\cat{M}.				\\
\end{eqnarray*}

\item A monoidal category \cat{C} with only one object is\ldots\ well, if we
forget the monoidal structure for a moment then we've already seen that
it's a monoid whose elements are the morphisms of \cat{C} and whose
multiplication is the composition in \cat{C}.  But the monoidal structure on
\cat{C} provides not only a tensor product for objects, but also a tensor
product for morphisms: so the set of morphisms of \cat{C} has a second
multiplication on it, $\otimes$.  So a one-object monoidal category is a set
$M$ equipped with two monoid structures which are in some sense compatible
(because of the axioms on a monoidal category).  And there's a well-known
result (the Eckmann-Hilton argument) saying that in this situation, the two
multiplications are in fact equal and commutative.  So: a one-object monoidal
category is a commutative monoid.  

This is essentially the same argument often used to prove that the higher
homotopy groups are abelian, or that the fundamental group of a topological
group is abelian.  In fact, we can deduce that $\pi_2$ is abelian from our
`results' so far:

\textbf{Corollary:} $\pi_2(X,x_0)$ is abelian (for a space $X$ with basepoint
$x_0$).  For the 2-category $\Pi_2(X)$ has a sub-2-category whose only 0-cell
is $x_0$, whose only 1-cell is the constant path at $x_0$, and whose 2-cells
are all the possible ones from $\Pi_2(X)$---that is, are the homotopies from
the constant path to itself, that is, are the elements of $\pi_2(X,x_0)$.
This sub-2-category is a 2-category with only one 0-cell and one 1-cell,
i.e.\ a monoidal category with only one object, i.e.\ a commutative monoid;
and this monoid is exactly $\pi_2(X,x_0)$.

\item Next consider a 3-category with only one 0-cell and one 1-cell.  We
haven't looked at (weak) 3-categories in enough detail to work this out
properly, but it turns out that such a 3-category is the same thing as a
braided monoidal category.  By definition, a \demph{braided monoidal
category} is a monoidal category equipped with a map (a \demph{braiding})
\[
A\otimes B \goby{\beta_{A,B}} B\otimes A
\]
for each pair $A,B$ of objects, satisfying axioms \emph{not} including that 
\[
(A\otimes B \goby{\beta_{A,B}} B\otimes A \goby{\beta_{B,A}} A\otimes B)
= 1.
\]
The canonical example of a braided monoidal category (in fact, the braided
monoidal category freely generated by a single object) is \fcat{Braid}.  This
has:
\begin{itemize}
\item objects: natural numbers $0, 1, \ldots$
\item morphisms: braids, e.g.\ 
\[
\begin{array}[c]{c}
\piccy{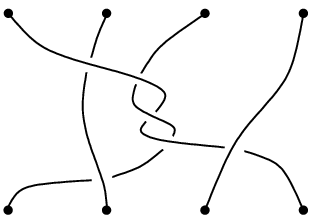}
\end{array}
\diagspace\diagspace\diagspace
\begin{diagram}[height=1cm]
4 \\ \dTo \\ 4 \\
\end{diagram}
\]
(taken up to deformation); there are no morphisms $m \go n$ when $m \neq n$
\item tensor: placing side-by-side (which on objects means addition)
\item braiding: left over right, e.g.
\[
\begin{array}[c]{c}
\piccy{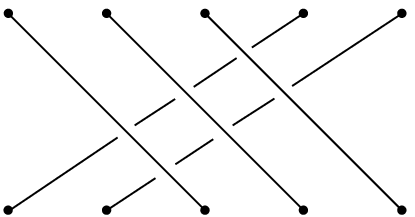}
\end{array}
\diagspace\diagspace
\begin{diagram}[height=1cm]
3+2 \\ \dTo<{\beta_{3,2}} \\ 2+3 \\
\end{diagram}
\]
(Notice how $\beta_{n,m} \of \beta_{m,n}$ is not the identity braid.)
\end{itemize}

\item We're getting way out of our depth here, but nevertheless: we've
already considered categories ($=$ categories), 2-categories which are only
interesting in the top two dimensions ($=$ monoidal categories), and
3-categories which are only interesting in the top two dimensions ($=$
braided monoidal categories).  What next?  For $r\geq 4$, an $r$-category
with only one $i$-cell for each $i<r-1$ is (people think) the same as a
symmetric monoidal category (i.e.\ a braided monoidal category in which
$\beta_{B,A} \of \beta_{A,B} = 1$ for all $A,B$).  So the situation's
stabilized\ldots\ and this is meant to make you start thinking of
stabilization phenomena in homotopy.
\end{itemize}

\paragraph*{The Big Picture}  Let's try to assemble this information on
degeneracies in a systematic way.  Define a \demph{$k$-monoidal $n$-category}
to be a $(k+n)$-category with only one $i$-cell for each $i<k$.  (It's clear
that this is going to be some kind of $n$-dimensional structure, as there are
only interesting cells in the top $(n+1)$ dimensions.)  Here's what
$k$-monoidal $n$-categories are for some low values of $k$ and $n$, laid out
in the so-called `periodic table'.  Explanation follows.

\begin{trivlist} \item
\hspace*{-10mm}
\setlength{\unitlength}{1.03mm}
\begin{picture}(128,88)(-1,-11)
\place{69}{76}{$n$}
\place{-1}{32}{$k$}
\put(2,66){\line(1,0){125}}
\put(10,0){\line(0,1){72}}
\place{6}{62}{$0$}
\place{6}{54}{$1$}
\place{6}{44}{$2$}
\place{6}{34}{$3$}
\place{6}{24}{$4$}
\place{6}{14}{$5$}
\place{6}{4}{$6$}
\place{24}{70}{$0$}
\place{24}{62}{set}
\place{24}{54}{monoid}
\place{24}{46}{commutative}
\place{24}{42}{monoid}
\place{24}{34}{\ditto}
\place{24}{24}{\ditto}
\place{24}{14}{\ditto}
\place{24}{4}{\ditto}
\place{54}{70}{$1$}
\place{54}{62}{category}
\place{54}{56}{monoidal}
\place{54}{52}{category}
\place{54}{46}{braided}
\place{54}{42}{mon cat}
\place{54}{36}{symmetric}
\place{54}{32}{mon cat}
\place{54}{24}{\ditto}
\place{54}{14}{\ditto}
\place{54}{4}{\ditto}
\place{84}{70}{$2$}
\place{84}{62}{2-category}
\place{84}{56}{monoidal}
\place{84}{52}{2-category}
\place{84}{46}{braided}
\place{84}{42}{mon 2-cat}
\place{84}{34}{X}
\place{84}{26}{symmetric}
\place{84}{22}{mon 2-cat}
\place{84}{14}{\ditto}
\place{84}{4}{\ditto}
\place{114}{70}{$3$}
\place{114}{62}{3-category}
\place{114}{56}{monoidal}
\place{114}{52}{3-category}
\place{114}{46}{braided}
\place{114}{42}{mon 3-cat}
\place{114}{34}{X}
\place{114}{24}{X}
\place{114}{16}{symmetric}
\place{114}{12}{mon 3-cat}
\place{114}{4}{\ditto}
\qbezier[10](24,54)(39,58)(54,62)
\qbezier[20](24,44)(54,53)(84,62)
\qbezier[30](24,34)(69,48)(114,62)
\qbezier[30](24,24)(69,39)(114,54)
\qbezier[30](24,14)(69,29)(114,44)
\qbezier[30](24,4)(69,19)(114,34)
\qbezier[20](54,4)(84,14)(114,24)
\qbezier[10](84,4)(99,9)(114,14)
\qbezier(38,-11)(45.5,-8.5)(53,-6)
\put(38,-11){\vector(-4,-1){0}}
\place{76}{-9}{:\ \ \ take just one-object things}
\end{picture}
\end{trivlist}

\emph{Running commentary:} in the first row ($k=0$), a $0$-monoidal
$n$-category is simply an $n$-category (it's not monoidal at all).

In the next row ($k=1$), a $1$-monoidal $n$-category is a monoidal
$n$-category, i.e.\ an $n$-category equipped with a tensor product which is
associative and unital up to equivalence of a suitable kind.  For instance, a
$1$-monoidal $0$-category is a one-object ($1$-)category, i.e.\ a monoid; and
a $1$-monoidal $1$-category is a one-object $2$-category, i.e.\ a monoidal
category.  We didn't look at the case of one-object $3$-categories, but they
turn out to be monoidal $2$-categories.  (What's a monoidal 2-category?
Well, it's a one-object 3-category\ldots\ or a direct definition can be
supplied.)  We see from these examples, or the general definition of
$k$-monoidal $n$-category, that going in the direction $\swarrow$ means
restricting to the one-object structures.

Now look at the third row ($k=2$): we've seen that a degenerate monoidal
category is a commutative monoid and a doubly-degenerate 3-category is a
braided monoidal category.  It's customary to keep writing `braided monoidal
$n$-category' all along the row, but you can regard this as nothing more than
name-calling.

Next consider the first \emph{column} ($n=0$).  A one-object braided monoidal
category is going to be a commutative monoid with a little extra data (for
the braiding) obeying some axioms, but it turns out that we don't actually
get anything new: in some sense, `you can't get better than a commutative
monoid'.  This gives the entry for $k=3, n=0$, and you get the same thing all
the way down the rest of the column.  

A similar story applies in the second column ($n=1$).  We saw---or rather, I
claimed---that for $k\geq 3$, a $k$-monoidal $1$-category is just a symmetric
monoidal category.  So again the column stabilizes, and again the point of
stabilization is `the most symmetric thing possible'.

The same goes in subsequent columns.  The X's could be replaced by more
terminology, e.g.\ the first is sometimes called `sylleptic monoidal
2-category', but it doesn't matter what that means.

The main point is that the table stabilizes for $k\geq n+2$---just like
$\pi_{k+n}(S^k)$.  So if you overlaid a table of the homotopy groups of
spheres onto the table above then they'd stabilize at the same points.  There
are arguments to see why this should be so (and I remind you that this is all
very informal and by no means completely understood).  Roughly, the fact that
the prototypical braided monoidal category \fcat{Braid} is not symmetric
comes down to the fact that you can't usually translate two 1-dimensional
affine subspaces of 3-dimensional space past each other; and this is the same
kind of dimensional calculation as you make when proving that the homotopy
groups of spheres stabilize.

\paragraph{Answer} to the initial question\ldots

\begin{description}
\item[A] Every weak 2-category is equivalent to a strict one.
In particular, every (weak) monoidal category is equivalent to a strict one.
So, for instance, we can pretend that the monoidal category of abelian groups
is strict, and so that $\otimes$ is strictly associative.

\item[B] \emph{Not} every weak 3-category is equivalent to a strict one.  In
fact we've already seen a counterexample.  For
\begin{itemize}
\item a weak 3-category with one $0$-cell and one $1$-cell is a braided
monoidal category
\item a strict 3-category with one $0$-cell and one $1$-cell is a strict
symmetric monoidal category
\item if a braided monoidal category is equivalent to a symmetric monoidal
category then it's symmetric.
\end{itemize}
So any braided monoidal category which is not symmetric is a weak 3-category
not equivalent to a strict one.  An example (the canonical example?) of this
is \fcat{Braid} itself; and the fact that \fcat{Braid} is not symmetric says
exactly that the overpass 
$\begin{array}[c]{c}
\piccy{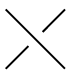} 
\end{array}$
can't be deformed to
the underpass 
$\begin{array}[c]{c}
\piccy{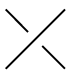} 
\end{array}$
in $\Bbb{R}^3$.
\end{description}

\pagebreak

\section{Further Reading}

There are many other things written on higher-dimensional category theory.
I'll only mention a few.

Another introduction to the theory, with many similar themes to this one, is
the first half of Baez's \cite{BzIN} paper.  Baez and collaborators have
written other interesting things on the interaction (still to be made
rigorous, mostly) between topology and higher-dimensional category theory,
and in particular his \cite{BDHDA0} paper with Dolan contains more on TQFT,
the periodic table, and stabilization.  

Specifically 2-categorical approaches to TQFT can be found in Tillmann
\cite{Till} and Khovanov \cite{Khov}.  $n$-vector spaces are to be found in
Kapranov and Voevodsky \cite{KV}, and their possible role in topological
field theory is discussed in Lawrence \cite{Law}.

Grothendieck puts the case that tame topology is really the study of
$\infty$-groupoids in his epic \cite{Gro} letter to Quillen.

One interesting idea that I didn't mention in this talk is that $n$th
cohomology should have coefficients in an $n$-category (as opposed to an
abelian group).  This is explained in the Introduction to Street's paper of
\cite{StAOS}.

A serious and, of course, highly recommended survey of the proposed
definitions of weak $n$-category, including ten such definitions, is my own
\cite{SDN} paper.  Some other occupants of the zoo of higher-dimensional
structures, including generalized operads and multicategories, are to be
found in my \cite{OHDCT}, the last chapter of my \cite{SHDCT}, and the second
half of Baez's \cite{BzIN}.

\it
\noindent Department of Pure Mathematics, University of Cambridge\\
Email: \url{leinster@dpmms.cam.ac.uk}\\
Web: \url{http://www.dpmms.cam.ac.uk/$\sim$leinster}


\begin{thebibliography}{2000}		

\small

\bibitem[1997]{BzIN}
John C. Baez, 
An introduction to $n$-categories.  
In \emph{7th Conference on Category Theory and Computer Science}, ed. 
E. Moggi and G. Rosolini, Springer Lecture Notes in Computer Science vol. 
1290 (1997). 
E-print \url{q-alg/9705009}.

\bibitem[1995]{BDHDA0}
John C. Baez, James Dolan,
Higher-dimensional algebra and topological quantum field theory.
\emph{Journal of Mathematical Physics} 36 (11) (1995), 6073--6105.
E-print \url{q-alg/9503002}(without pictures).

\bibitem[1984]{Gro}
A. Grothendieck,
Pursuing stacks (1984).
Manuscript.

\bibitem[1994]{KV}
M. M. Kapranov, V. A. Voevodsky,
2-categories and Zamolodchikov tetrahedra equations (1994).
In \emph{Algebraic Groups and their Generalizations: Quantum and
Infinite-Dimensional Methods}, Proceedings of Symposia in Pure Mathematics 56,
Part 2, AMS, 177--260.

\bibitem[2001]{Khov}
Mikhail Khovanov,
A functor-valued invariant of tangles (2001).
E-print \url{math.QA/0103190}.

\bibitem[1996]{Law}
R. J. Lawrence,
An introduction to topological field theory (1996).
In \emph{The Interface of Knots and Physics}, 
Proceedings of Symposia in Applied Mathematics 51, AMS, 89--128.
Also available via \url{http://www.ma.huji.ac.il/$\sim$ruthel}.

\bibitem[1998]{SHDCT}
Tom Leinster,
Structures in higher-dimensional category theory (1998).
Available via \url{http://www.dpmms.cam.ac.uk/$\sim$leinster}.

\bibitem[2000]{OHDCT}
Tom Leinster,
Operads in higher-dimensional category theory.
PhD thesis, University of Cambridge, 2000.  
E-print \url{math.CT/0011106}.

\bibitem[2001]{SDN}
Tom Leinster,
A survey of definitions of $n$-category (2001).
E-print \url{math.CT/01xxxxx}(to appear shortly).

\bibitem[1987]{StAOS}
Ross Street,
The algebra of oriented simplexes.
\emph{Journal of Pure and Applied Algebra} 49 (1987), no.~3, 283--335.

\bibitem[1998]{Till}
Ulrike Tillmann,
$\mathcal{S}$-structures for $k$-linear categories and the definition of a
modular functor.
\emph{Journal of the London Mathematical Society (2)} 58 (1998), no.~1,
208--228.
E-print \url{math.GT/9802089}.

\end{thebibliography}
\end{document}